\documentclass[12pt, a4paper]{amsart}
\usepackage{amsmath}
\usepackage{geometry,amsthm,graphics,tabularx,amssymb,shapepar}
\usepackage{latexsym,amsfonts,amssymb,amsmath,amsxtra,bbm,color,mathrsfs}
\usepackage{amscd}
\usepackage{hyperref}
\usepackage{pdfsync}
\usepackage[all]{xypic}

\newcommand{\K}{\mathsf{k}}

\newcommand{\BC}{{\mathbb {C}}}

\newcommand{\BP}{{\mathbb {P}}}

\newcommand{\BR}{{\mathbb {R}}}

\newcommand{\CF}{{\mathcal {F}}}

\newcommand{\CO}{{\mathcal {O}}}

\newcommand{\CS}{{\mathcal {S}}}

\newcommand{\RH}{{\mathrm {H}}}

\newcommand{\RN}{{\mathrm {N}}}

\newcommand{\RT}{{\mathrm {T}}}

\newcommand{\GL}{{\mathrm{GL}}}

\newcommand{\Hom}{{\mathrm{Hom}}}

\renewcommand{\Re}{{\mathrm{Re}}}

\newcommand{\diag}{\operatorname{diag}}

\newcommand{\oL}{\operatorname{L}}

\newcommand{\oO}{\operatorname{O}}

\newcommand{\oU}{\operatorname{U}}
\newcommand{\oZ}{\operatorname{Z}}

\newcommand{\C}{\mathbb{C}}
\newcommand{\R}{\mathbb R}

\newcommand{\abs}[1]{\lvert#1\rvert}

\newcommand{\la}{\langle}

\newcommand{\be}{\begin {equation}}
\newcommand{\ee}{\end {equation}}
\newcommand{\bee}{\begin {equation*}}
\newcommand{\eee}{\end {equation*}}

\theoremstyle{Theorem}

\theoremstyle{Theorem}

\theoremstyle{Theorem}
\newtheorem{lem}{Lemma}[section]

\newtheorem{prpl}[lem]{Proposition}

\theoremstyle{Theorem}

\theoremstyle{Plain}

\theoremstyle{remark}

\theoremstyle{remark}
\newtheorem*{rremark}{Remark}

\theoremstyle{Definition}
\newtheorem{dfn}{Definition}[section]

\newtheorem{prpd}[dfn]{Proposition}
\newtheorem{thmd}[dfn]{Theorem}

\newtheorem{exampled}[dfn]{Example}

\newcommand{\of}[1]{\left(#1\right)}
\def\dd{\mathrm{d}}
\def\la{\lambda}
\def\La{\Lambda}
\def\ga{\gamma}
\def\si{\sigma}
\def\om{\omega}
\def\vep{\varepsilon}

\begin{document}

\title[Rankin-Selberg integrals]{Rankin-Selberg integrals for principal series representations of $\GL(n)$}

\author[D. Liu]{Dongwen Liu}
\address{School of Mathematical Sciences, Zhejiang University, Hangzhou, Zhejiang, P. R. China}
\email{maliu@zju.edu.cn}

\author[F. Su]{Feng Su}
\address{Department of Pure Mathematics, Xi’an Jiaotong–Liverpool University, Suzhou 215123, China}
\email{feng.su@xjtlu.edu.cn}

\author[B. Sun]{Binyong Sun}
\address{Institute for Advanced Study in Mathematics, Zhejiang University, Hangzhou, Zhejiang, P. R. China}
\email{sunbinyong@zju.edu.cn}




\begin{abstract}
We prove that the local Rankin--Selberg integrals for principal series representations of the general linear groups agree with certain simple integrals over the Rankin--Selberg subgroups, up to certain constants given by the local gamma factors.

\end{abstract}

\maketitle

\section{Introduction and the main results}

Rankin--Selberg integrals provide a powerful tool in the study of automorphic representations and L-functions. Explicit calculations of the local
Rankin--Selberg integrals are often desirable for arithmetic applications.  
The goal of this note is to show that the local Rankin--Selberg integrals for principal series representations of the general linear groups agree with certain simple integrals over the Rankin--Selberg subgroups, up to certain explicit constants given by the local gamma factors.

Fix an arbitrary local field $\K$. Let $G:=\GL_n(\K)$ ($n\geq 2$). Let $B=AN$ be the Borel subgroups of $G$ of the upper-triangular matrices, where $N$ is the unipotent radical of $B$, and $A$ is the subgroup of the diagonal matrices. Similarly, let $\bar B=A\bar N$ be the Borel subgroup of $G$ of the lower-triangluar matrices, where $\bar N$ is the unipotent radical of $\bar B$.

The first Rankin-Selberg subgroup of $G$ is the group $R$ consisting of all matrices of the form 
\[
  \begin{bmatrix} a & u \\
0 & h \end{bmatrix}\in G
\]
such that $ a\in \K^\times$, $h$ is upper-triangular unipotent, and $u$ is a row vector whose first entry equals $0$. We put 
\[N':=R\cap N\qquad \textrm{and}\qquad A':=R\cap A\cong \K^\times \]
so that $R=A'N'$.

Fix a non-trivial unitary character $\psi: \K\rightarrow \C^\times$, and equip $\K$ with the  self-dual Haar measure $\mathrm d x$ associated to $\psi$.  Write $\abs{\,\cdot \,}_\K$ for the normalized absolute value on $\K$. We equip the following Haar measures on $\K^\times$, $N$ and $R$ respectively:
\begin{itemize}
 \item $\mathrm d^\times a:=\abs{a}_\K^{-1} \mathrm d a$, $\quad a\in \K^\times$; \smallskip
\item
     $\mathrm d u:=\prod_{1\leq i<j\leq n} \mathrm d u_{i,j}$, $\quad u=[u_{i,j}]_{1\leq i,j\leq n} \in N$;\smallskip
\item
      $\mathrm d_{\mathrm r} g:=\prod_{1\leq i<j\leq n, j\neq 2} \mathrm d g_{i,j}\cdot \mathrm d^\times g_{1,1}$, $\quad g=[g_{i,j}]_{1\leq i,j\leq n} \in R$. 
\end{itemize} 

Using $\psi$ we define the following  character of $N$:
\[
\psi_N: N\rightarrow \C^\times, \qquad [u_{i,j}]_{1\leq i,j\leq n}\mapsto \psi\of{\sum_{i=1}^{n-1}u_{i,i+1}}.
\]
Write $\Hom(A,\C^\times)$ for the set of all characters of $A$, which is a complex Lie group of dimension $n$. Let $\sigma=\sigma_1\otimes \sigma_2\otimes \cdots \otimes \sigma_n\in \Hom(A,\C^\times)$, where $\sigma_1, \sigma_2, \cdots, \sigma_n$ are characters of $\K^\times$ . View $\sigma$ as a character of $\bar B$ through the trivial extension to $\bar N$, and define the principal series representation 
\[
  I(\sigma):=\mathrm{Ind}_{\bar B}^G \sigma.
\] 
Recall that $I(\sigma)$ consists of smooth functions $f: G\to\BC$ such that 
\[
f(\bar{b}\cdot g)=\sigma(\bar{b})\cdot \bar{\rho}(\bar{b}) \cdot f(g),\quad \textrm{for all } \, b\in \bar{B}, \ g\in G,
\]
where
\[
\bar{\rho}=\abs{\,\cdot\,}_\K^{\frac{1-n}{2}}\otimes \abs{\,\cdot\,}_\K^{\frac{3-n}{2}}\otimes \ldots \otimes \abs{\,\cdot\, }_\K^{\frac{n-1}{2}}\in \Hom(A,\C^\times),
\]
and that $G$ acts on $I(\sigma)$ through the right translations. 
In the non-archimedean case,  $ I(\sigma)$ is equipped with the finest locally convex topology such that all seminorms on it are continuous. In the archimedean case,  $ I(\sigma)$ is a Fr\'echet space under the smooth topology.

 It is known that  there is a unique element $\lambda:=\lambda_\sigma\in \Hom_N(I(\sigma), \psi_N)$ such that
\[
  \lambda(f)=\int_N f(u) \psi_N^{-1}(u)\,\mathrm d u
\]
for all $f\in I(\sigma)$ such that $f|_N\in \CS(N)$  (see \cite[Theorem 15.4.1]{Wa92}). Here and henceforth, we write $\CS(X)$ for the space of Schwartz functions on $X$ when $X$ is a Nash manifold (see \cite{AG08}), and the space of compactly supported locally constant functions on $X$ when $X$ is a totally disconnected locally compact topological space. All functions are complex valued unless otherwise specified.

For every $a\in \K^\times$, write
\[
[a]:=\diag(a,1, \cdots,1)\in G.
\]
Then $A'=\{[a] \ | \ a\in \K^\times\}$. 
For every $f\in I(\sigma)$ and $s\in \C$, the local Rankin--Selberg integral is define to be 
\[
  \oZ_s(f):=\int_{\K^\times}\la\of{[a].f}\abs{a}_\K^{s-\frac{n-1}{2}}\dd^\times a.
\]
For every $s\in\C$, define a character 
\[
\psi_s:\,R\to\C^\times,\quad u'\cdot[a]\mapsto\psi_N(u')\cdot \abs{a}_\K^{\frac{n-1}{2}-s},\quad u'\in N',\ a\in\K^\times.
\]
For every  character $\om$ of $\K^\times$, let $\oL(s,\om)$ denote the local $\oL$-function of $\om$. Write \[
\mathrm L(s, \sigma):=\prod_{i=1}^n \mathrm L(s, \sigma_i),
\]
 which is a meromorphic function on $\C$. Note that  $\frac{1}{\mathrm L(s, \sigma)}$ is an entire function. 

Some basic properties of the Rankin--Selberg integrals are summarized in the following proposition.  See \cite[Section 5.3]{J09} and \cite[Section 8.3]{JPSS83}.

\begin{prpd} \label{prop1.1}
There is a real number $C_\sigma$ with the following properties.
\begin{itemize}
\item 
For all $s\in \C$ with the real part  $\Re(s)>C_\sigma$,  the integral $\oZ_s(f)$ converges absolutely for all $f\in I(\sigma)$. 
\item The map 
\[
\{s\in \C\mid\Re(s)>C_\sigma\} \times I(\sigma)\to \C,\qquad  (s, f)\mapsto \oZ_s(f),
\]
  extends to the multiplication of $\oL(s,\sigma)$ with a continuous map 
\[
  \oZ^\circ: \C\times I(\sigma) \to \C
  \]
that is holomorphic on the first variable and linear on the second variable. Moreover, for every $s\in \C$,
\[
\oZ^\circ(s, \cdot) \in \Hom_R(I(\sigma), \psi_s).
\]
 \end{itemize}

\end{prpd}

On the other hand, we set
\[w_1=\diag\left(\begin{bmatrix}1&1\\0 &1\end{bmatrix}, 1, \cdots, 1\right)\in G,\]
and  define the integral 
\[
\La_s(f):= \int_Rf(w_1g)\psi_s^{-1}(g)\dd_r g,\qquad s\in \C,\  f\in I(\sigma). 
\]

Fore each $i=1,2, \cdots, n$, write $\nu_i$ for the real number such that $\abs{\sigma_i}=\abs{\,\cdot\,}_\K^{\nu_i}$ as positive characters of $\K^\times$. 
Put
\[
 \Omega_\sigma:=\{s\in \C\mid -\nu_2<\Re(s) <1-\nu_1\}. 
\]

\begin{thmd}\label{thm1}
Assume that 
\be \label{nu}
\left\{\begin{aligned}
 & \max\{\nu_1, \nu_2\}<\nu_3<\cdots<\nu_{n-1}<\nu_n,\\
&  \nu_1<\nu_2+1.
\end{aligned}
\right.
\ee
Then for all $s\in \Omega_\sigma$, the integral $\La_s(f)$ converges absolutely for all $f\in I(\sigma)$. Moreover, the map
\[
  \Omega_\sigma\times  I(\sigma) \rightarrow \BC, \qquad (s,f)\mapsto \La_s(f)
\]
is continuous,    holomorphic on the first variable, and linear on the second variable.  
\end{thmd}

Under the assumptions of Theorem \ref{thm1}, we get an element $\La_s\in \Hom_R(I(\sigma), \psi_s)$ for every $s\in \Omega_\sigma$. 
It is clear that $\La_s\neq 0$. 

\begin{thmd}\label{unique} For all but countably many $s\in \BC$, 
\[
\dim \Hom_R(I(\sigma),\psi_s) = 1.
\]

\end{thmd}

In the  non-archimedean case,  Theorem \ref{unique} is proved in \cite[Proposition 2.11]{JPSS83} in a more general setting.  We will prove Theorem \ref{unique} in the archimedean case by using the theory of Schwartz homologies in \cite{CS21}. 

For  every $s\in \Omega_\sigma$, Theorem \ref{unique} implies  that $\La_s$  equals $\oZ_s$ up to a scalar multiplication. 
More precisely, we will prove the following result.

\begin{thmd}\label{thm2}
Let the notation and assumptions be as in Theorem \ref{thm1}. Then
\be \label{maineq}
\La_s(f)=\ga(s, \sigma_1, \psi)\cdot\oZ_s(f)
\ee
 for all $f\in I(\sigma)$ and $s\in \Omega_\sigma$.
 \end{thmd}
 
 \begin{rremark}
  Here $\ga(s, \sigma_1, \psi)$ denotes the usual local gamma factor which will be  recalled later. The right hand side of \eqref{maineq} is holomorphic in $s\in \Omega_\sigma$, which can be seen from the equalities
\[
 \begin{aligned}
  \gamma(s, \sigma_1, \psi)\cdot \oZ_s(f)  
   & =   \varepsilon(s, \sigma_1, \psi) \cdot \frac{\oL(1-s, \sigma_1^{-1})}{\oL(s, \sigma_1)} \cdot \oL(s, \sigma) \cdot \oZ^\circ (s, f) \\
   & =  \varepsilon(s, \sigma_1, \psi) \cdot \oL(1-s, \sigma_1^{-1}) \cdot \prod^n_{i=2} \oL(s, \sigma_i) \cdot \oZ^\circ (s, f).
 \end{aligned}
 \]
 Here $ \oZ^\circ $ is as in Proposition \ref{prop1.1}.
 \end{rremark}

We recall the definition of local gamma factor following
\cite{T79, J79, K03}.  Given a character $\omega$ of $\K^\times$, the Tate's local zeta integral (\cite{T50}) is defined by
\[
\oZ(s,\om, f)=\int_{\K^\times}f(x)\om(x)\abs{x}_\K^s\, \dd^\times x,\quad  f \in\CS(\K),\]
which converges absolutely when $\Re(s)>-{\rm ex}(\omega)$. Here ${\rm ex}(\omega)$ is the real number such that $|\omega|=\abs{\,\cdot\,}_\K^{{\rm ex}(\omega)}$ as 
positive characters of $\K^\times$.  The local epsilon factor attached to $\om$ and $\psi$ will be denoted by $\vep(s,\om,\psi)$, which is defined by  the local functional equation (\cite{T50})
\begin{equation}\label{fe}
\frac{\oZ(1-s,\om^{-1}, \CF_\psi(f))}{\oL(1-s,\om^{-1})}=\vep(s,\om,\psi)\cdot\frac{\oZ(s,\om, f)}{\oL(s,\om)},\quad  f \in \CS(\K),
\end{equation}
where $\CF_\psi(f)\in \CS(\K)$ is the Fourier transform of $f$ with respect to $\psi$ defined by
\[
\CF_\psi(f)(x):=\int_{\K}f(y)\psi(xy)\dd y,\quad x\in \K.
\]
The meromorphic function 
\[
\ga(s,\om,\psi):=\vep(s,\om,\psi)\cdot\frac{\oL(1-s,\om^{-1})}{\oL(s,\om)}
\] 
is called the local gamma factor attached to $\om$ and $\psi$.

The functional $\La_s$  takes  a simpler form than $\oZ_s$, and we expect that Theorem \ref{thm2} will be useful for the study of global L-functions. 


\section{Proof of Theorem \ref{thm1}}

Let $f\in I(\sigma)$. We first  rewrite $\La_s(f)$  formally as an integral over $N$. Let $N_\alpha$ be the root subgroup of $N$ corresponding to the positive simple root $\alpha:=e_1-e_2$, that is, $N_\alpha$ consists of the matrices of the form
\[
\mathrm u(x):={\rm diag}\left(\begin{bmatrix} 1 & x \\ 0 & 1\end{bmatrix}, 1, \cdots, 1\right),\quad x\in \K.
\]
Then $N=N_\alpha N'=N'N_\alpha$. In view of the equalities
\[
[a]^{-1}w_1[a]=[a]^{-1}\mathrm u(1)[a]=\mathrm u(a^{-1}),\quad  a\in \K^\times,
\]
 we have that
\begin{eqnarray*}
\La_s(f)&=&\int_{\K^\times}\int_{N'}f(w_1u'[a])\psi_s^{-1}(u'[a])\dd u'\dd^\times a\\
&=&\int_{\K^\times}\int_{N'}f([a]\cdot \mathrm u(a^{-1})\cdot[a]^{-1}u'[a])\psi_s^{-1}(u'[a])\dd u'\dd^\times a\\
&=&\int_{\K^\times}\int_{N'}\si_1(a)|a|_{\K}^{-\frac{n-1}{2}} f(\mathrm u(a^{-1})\cdot[a]^{-1}u'[a])\psi_s^{-1}(u'[a])\dd u'\dd^\times a.
\end{eqnarray*}
The change of variable $u'\mapsto [a]u'[a]^{-1}$ does not affect the value of $\psi_s(u'[a])$, and we have that
\[
\big |\det \left(\left.\textrm{Ad}([a])\right|_{{\rm Lie  \ }N'}\right)\big |_\K=|a|_\K^{n-2}.
\] 
Hence 
\begin{equation} \label{eq2}
\begin{aligned}
\La_s(f)&=\int_{\K^\times}\int_{N'}\si_1(a)|a|_{\K}^{\frac{n-3}{2}} f(\mathrm u(a^{-1})\cdot u')\psi_s^{-1}(u'[a])\dd u'\dd^\times a\\
&= \int_{\K^\times}\int_{N'}\si_1(a)|a|_{\K}^{s-1} f(\mathrm u(a^{-1})\cdot u')\psi_N^{-1}(u')\dd u'\dd^\times a\\
&= \int_{\K^\times}\int_{N'}\si_1^{-1}(a)|a|_{\K}^{1-s} f(\mathrm u(a)\cdot u')\psi_N^{-1}(u')\dd u'\dd^\times a,
\end{aligned}
\end{equation}
where we made the change of variable $a\mapsto a^{-1}$ in the last equality.

Introduce an open subset of $N$,
\begin{equation} \label{N*}
N^*:=\{u\in [u_{i, j}]\in N \ | \ u_{1,2}\neq 0\},
\end{equation} 
and for  $s\in \C$ define a function $h_s$ on $N^*$ by
\[
h_s(u):=|u_{1,2}|_\K^{-\nu_1-s},\quad u=[u_{i,j}]\in N^*.
\]
Let $K$ be the standard maximal compact subgroup of $G$, namely
\[
 K:=\begin{cases}
 \oO(n), \quad &\textrm{if }\K\cong \R;\\
  \oU(n), \quad &\textrm{if }\K\cong \C;\\
   \GL_n(\mathfrak o_\K), \quad &\textrm{if $\K$ is nonarchimedean,}
 \end{cases}
\]
where $\mathfrak o_\K$ denotes the ring of integers in $\K$. 
Put
\[
\|f\|_K:=\max_{k\in K}|f(k)|.
\]
Then for $a\in \K^\times$ and $u' \in N'$ we have that
\be \label{bound}
\left| \si_1^{-1}(a)|a|_{\K}^{-s} f(\mathrm u(a)\cdot u')\psi_N^{-1}(u')\right| \leq \|f \|_K\cdot \varphi_\nu(\mathrm u(a)\cdot u') \cdot h_{{\rm Re}(s)}(\mathrm u(a)\cdot u'),
\ee
where $\nu:=\abs{\, \cdot \, }_\K^{\nu_1}\otimes\cdots\otimes \abs{\, \cdot \, }_\K^{\nu_n} \in \Hom(A,\C^\times)$, and  $\varphi_\nu$ is the spherical vector of $I(\nu)$ satisfying that $\varphi_\nu|_K\equiv 1$.
Note that $\varphi_\nu$ is positive valued. Define the integral 
\[
\eta_\nu(s):=\int_N \varphi_\nu(u) h_s(u) \dd u.
\]


To prove Theorem \ref{thm1}, we first prove the following result.

\begin{prpl} \label{prop2.1}
Assume that \eqref{nu} holds.
Then the integral $\eta_\nu(s)$ converges absolutely and uniformly on compact subsets of $\Omega_\sigma$.
\end{prpl}

\begin{proof}
Let $G_\alpha=\GL_2$ embedded into the top-left corner of $G$, with Iwasawa decomposition 
\[
G_\alpha=\bar{N}_\alpha A_\alpha K_\alpha,
\]
where $\bar{N}_\alpha$ is  the lower-triangular unipotent subgroup, $A_\alpha$ is the diagonal torus and $K_\alpha$ is the standard maximal compact subgroup of $G_\alpha$. Write the Iwasawa decomposition 
of $\mathrm u(x)$, $x\in \K$, accordingly as
\[
\mathrm u(x)=\bar{n}(x) \cdot a(x) \cdot k(x).
\]
Let $W$ be the subgroup of permutation matrices in $G$, and $s_\alpha\in W$ be the simple reflection corresponding to $\alpha$. Since $N_\alpha$ and $s_\alpha$  normalize $N'$, we see that $\bar{N}_\alpha={\rm Ad}(s_\alpha)(N_\alpha)$ normalizes $N'$ as well. It is clear that for $x\in \K$ the map
\[
{\rm Ad}(\bar{n}(x)): N'\longrightarrow N'
\]
is unimodular. Let $\Phi(N', A)$ be the set  of  roots of $A$ in $N'$. The roots in $\Phi(N',A)$ add up to
\[
\mu:=2\rho-\alpha,
\]
where $\rho$ acts on $A$ by the character $\abs{\,\cdot\,}_\K^{\frac{n-1}{2}}\otimes  \ldots \otimes \abs{\,\cdot\, }_\K^{\frac{1-n}{2}}$. Then we have that
\[
\big | \det({\rm Ad}(a(x))|_{{\rm Lie} \ N'})\big |_\K= |\mu(a(x))|_\K.
\]
 It follows  that we have formally 
\be \label{eta}
\begin{aligned}
\eta_\nu(s) &=  \int_{\K} \int_{N'} \varphi_\nu(u' \cdot \mathrm u(x))h_s(\mathrm u(x))\dd u' \dd x \\
&=  \int_{\K} (\nu\cdot \bar{\rho}\cdot \mu)(a(x)) |x|_\K^{-\nu_1-s} \dd x \cdot \int_{ N'} \varphi_\nu(u')\dd u'. 
\end{aligned}
\ee

Let $\Phi$ be the standard root system  of $\GL_n$. For $w\in W$, let $N_w$ be the subgroup of $N$ corresponding to the set of positive roots
\begin{equation} \label{phiw}
\Phi_w:= \Phi_+\cap w^{-1}\Phi_-.
\end{equation}
Then it is clear that $N'=N_{w_0s_\alpha}$, where $w_0\in W$ is the anti-diagonal matrix, and that $\Phi_{w_0s_\alpha}=\Phi(N', A)$.
In view of this, the second integral in the last product in \eqref{eta} is among the Harish-Chandra's $c$-functions
(see e.g. \cite{L67, GGPS90} and  \cite[Chapter IV]{H84}), and it converges if and only if 
\[
 \langle \nu, \beta^\vee\rangle <0
\]
for all $\beta\in \Phi(N', A)$, where $\beta^\vee$ denotes the coroot of $\beta$, and $\langle, \rangle$ is the natural pairing between $X^*(A)_{\BC}$ and $X_*(A)_\BC$. The above condition reads
\[
\max\{\nu_1, \nu_2\}<\nu_3<\cdots<\nu_{n-1}<\nu_n.
\]
In this case, up to a positive scalar depending on the Haar measure,  the Gindikin-Karpelevich formula  (see \cite[Section 4]{L67}) gives that
\[ \int_{ N'} \varphi_\nu(u')\dd u'=\prod_{\beta\in \Phi(N', A)} \frac{\zeta_\K(-\langle \nu, \beta^\vee\rangle)}{\zeta_\K(1-\langle\nu, \beta^\vee\rangle)}.\]
Here $\zeta_\K(s)$ is the zeta function of $\K$ (see \cite[Section 3]{T79}) defined by
\[
\zeta_\K(s)=\begin{cases}
\pi^{-s/2}\Gamma(s/2) ,&\text{if $\K=\R$};\\[.5em]
 2(2\pi)^{-s}\Gamma(s),&\text{if $\K=\C$};\\[.5em]
\dfrac{1}{1-q_\K^{-s}}, &\text{if $\K$ is nonarchimedean},
\end{cases}
\]
with $\Gamma(s)$  the usual Gamma function, and  $q_\K$  the cardinality of the residue field of $\frak{o}_\K$ in the last case. 

By the standard computation of Iwasawa decomposition for $\GL_2(\K)$, the first integral in the last product in \eqref{eta} is 
\begin{eqnarray*}
&&  \int_{\K} (\nu\cdot \bar{\rho}\cdot \mu)(a(x)) |x|_\K^{-\nu_1-s} \dd x \\
&=& \begin{cases}
\int_\K(1+\abs{x}_\K^2)^{\frac{\nu_1-\nu_2-1}{2}}\abs{x}_\K^{-\nu_1-s}\dd x,&\text{if $\K=\R$};\\[.5em]
\int_\K(1+\abs{x}_\K)^{\nu_1-\nu_2-1}\abs{x}_\K^{-\nu_1-s}\dd x,&\text{if $\K=\C$};\\[.5em]
\int_{\abs{x}_\K\leq 1}\abs{x}_\K^{-\nu_1-s}\dd x+\int_{\abs{x}_\K>1}\abs{x}_\K^{-\nu_2-s-1}\dd x,&\text{if $\K$ is nonarchimedean}.
\end{cases}
\end{eqnarray*}
It is easy to check that the above integral  converges absolutely and uniformly on  compact subsets of $\Omega_\sigma$. 

This proves the Proposition. 
\end{proof}

By \eqref{eq2}, \eqref{bound} and Proposition \ref{prop2.1},  we see  that $\La_s(f)$ converges absolutely and uniformly 
on compact subsets of $\Omega_\sigma$.  By the Weierstrass Theorem in complex analysis (see \cite[Chap. 5, Theorem 1]{A78}),  $\La_s(f)$ is holomorphic on $\Omega_\sigma$. 

For $(s, f)$, $(s', f')\in \Omega_\sigma\times I(\sigma)$, we have that
\[
\begin{aligned}
|\Lambda_s(f)-\Lambda_{s'}(f')| & \leq |\Lambda_s(f-f')| + | \Lambda_s(f')-\Lambda_{s'}(f')| \\
& \leq \| f-f'\|_K\cdot \eta_\nu({\rm Re}(s)) +\|f'\|_K \cdot\int_N \varphi_\nu(u) \left| h_s(u) - h_{s'}(u) \right|  \dd u.
\end{aligned}
\]
It follows from the above inequalities and Proposition \ref{prop2.1} that the map
\[
\Omega_\sigma\times I(\sigma)\to \C, \quad (s, f)\mapsto \Lambda_s(f)
\]
is continuous. This finishes the proof of Theorem \ref{thm1}.

\section{$R$-orbits and proof of Theorem \ref{unique}}

In this section we prove Theorem \ref{unique} regarding the uniqueness of Rankin--Selberg periods.

\subsection{$R$-orbits on the flag variety} 
Let $M:=\bar{B}\backslash G$ be the flag variety of Borel subgroups of $G$ with base point $x_0:=\bar{B}\in M$.
For every $g\in G$, denote the $R$-orbit of $x_0\cdot g\in M$ by
\[
\mathcal{O}_g:=\{x_0 \cdot (g r)\mid r\in R\}\subset M.
\]
We shall describe the $R$-orbits in $M$, or equivalently the $\bar{B}$-$R$ double cosets. 
To this end we make some preparations.  For every $w\in W$ we have the opposite Bruhat cell
\[
C_w:= x_0 \cdot wB=x_0\cdot  w U_w,
\]
where $U_w:=N_{w_0 w}$ is the subgroup of $N$ corresponding to the  following set of positive roots (cf. \eqref{phiw})
\[
\Psi_w: =\Phi_{w_0w}=\Phi_+\cap w^{-1}\Phi_+.
\]
Then $C_1=x_0\cdot N$ is open dense in $M$, and we recall the  Bruhat decomposition 
\[
M=\bigsqcup_{w\in W}C_w.
\]
Introduce the following subsets of $W$:
\[
\begin{aligned}
& W_+:=\{w\in W \ | \ w\alpha>0\}, \\
& W_-:=\{w\in W \ | \ w\alpha<0\}.
\end{aligned}
\]

\begin{prpd} \label{prop3.1}
We have the following disjoint union of $R$-orbits
\[
M=\left(\bigsqcup_{w\in W}\CO_w\right) \bigsqcup\left(\bigsqcup_{w\in W_+} \CO_{ww_1} \right).
\]
Moreover it holds that
\begin{itemize}
\item[(i)] if $w\in W_-$ then $C_w=\CO_w$;
\item[(ii)] if $w\in W_+$ then $C_w=\CO_w\sqcup \CO_{ww_1}$ and $\CO_{ww_1}$ is open dense in $C_w$. 
\end{itemize}
\end{prpd}

\begin{proof}
By the Bruhat decomposition, it suffices to prove (i) and (ii). Put
\[
N_\alpha^*:=N_\alpha\cap N^*=N_\alpha\setminus\{1\}=\{\mathrm u(a) \ | \ a\in \K^\times\}.
\] 
If $w\in W_-$, then
\[
C_w=x_0\cdot wN_\alpha N'=x_0 \cdot wN'= x_0 \cdot wR=\CO_w.
\] 
If $w\in W_+$, then 
\[
C_w=x_0\cdot wN=x_0\cdot w N' \sqcup x_0\cdot w N_\alpha^* N'.
\]
We have $x_0\cdot wN'=\CO_w$, and it is easy to verify that
\[
x_0\cdot wN_\alpha^*N'=x_0\cdot wA'N_\alpha^* N'=x_0\cdot w A'w_1 A'N'=x_0\cdot ww_1R=\CO_{ww_1},
\]
which finishes the proof.
\end{proof}

\begin{exampled}
If $n=2$, then $R=\K^\times$ acts on $M = \BP^1(\K)$ with three orbits $\{0\}$, 
$\{\infty\}$ and $\K^\times$.
\end{exampled}

\subsection{Proof of Theorem \ref{unique}} Assume that $\K$ is archimedean. Then we have the following result, which in particular implies 
Theorem \ref{unique}.

\begin{thmd}
For all but countably many $s\in \BC$, there is a topological linear isomorphism 
\[
\RH_i^\CS(R; I(\sigma)\otimes\psi_s^{-1})\cong \left\{\begin{array}{ll}\BC, & {\rm if \ }i=0;\\
\{0\}, & {\rm if \ }i \neq 0.
\end{array}
\right.
\]
\end{thmd}

Here $\RH^\CS_i$ indicates the Schwartz homology studied in \cite{CS21}.

\begin{proof}
The flag variety $M$ is naturally a $G$-Nash manifold, and we have 
\[
I(\sigma)=\Gamma^\zeta(M, {\sf E})
\]
for a certain tempered $G$-vector bundle ${\sf E}$ of rank one over $M$. Here $\Gamma^\zeta(M, {\sf E})$ is the Fr\'echet space of Schwartz sections of $\sf{E}$ defined as in \cite[Section 6.1]{CS21}. 

Denote $U:=\CO_{w_1}$ the unique open $R$-orbit in $M$, and $Z:=M\setminus U$ its complement. Then 
\[
\Gamma^\zeta(U, \left.{\sf E}\right|_U)\cong \CS(R).
\]
For $z\in M$ we have its stabilizer in $R$ given by
\[
R_z=R\cap \tilde{z}^{-1}\bar{B}\tilde{z},
\]
where $\tilde{z}$ is an arbitrary representative of $z$ in $G$. Write 
\[
\RN^*_z:=\frac{\RT^*_z(M)}{\RT^*_z(R.z)}\otimes_\BR\BC \quad (\RT_z^*\textrm{ stands for the cotangent space})
\]
for the complexified conormal space, and $\delta_{R/R_z}:=(\delta_R)|_{R_z}\cdot \delta_{R_z}^{-1}: R_z\to \BC^\times$ where $\delta$ stands for the modular character. 

We claim that for all but countably many $s\in \BC$, the condition of \cite[Theorem 1.15]{CS21} holds, namely for all $z\in Z$ and all integers $k\geq 0$, the trivial representation of $R_z$ does not occur as a subquotient of
\begin{equation}\label{symp}
{\sf E}_z\otimes {\rm Sym}^k(\RN_z^*)\otimes \delta_{R/R_z}\otimes \psi_s^{-1},
\end{equation}
where ${\sf E}_z$ is the fibre of ${\sf E}$ at $z$. 
By Proposition \ref{prop3.1}, we may assume that $z=w$, $w\in W$ or $z=ww_1$, $w\in W_+\setminus\{1\}$. 

 If $z=w$, $w\in W$, then $A'\subset R_z$. 
 
 If $z=ww_1$, $w\in W_+\setminus\{1\}$, then $w$ maps at least one of the positive simple roots
$\alpha_i=e_i-e_{i+1}$, $i=2,\ldots, n-1$ to a negative root. At this point we consider two cases separately:
\begin{itemize}
\item $w\alpha_i<0$ for some $i\geq 3$. Then $w_1\in N_\alpha$ commutes with $N_{\alpha_i}$, hence 
$N_{\alpha_i}\subset R_z$;
\item $w\alpha_2<0$. Let $S$ be the subgroup of $R$ consisting of the matrices
\[
{\rm diag}\left(\begin{bmatrix} 1 & 0 & -x \\ 0 & 1 & x \\ 0 & 0 & 1\end{bmatrix}, 1, \cdots, 1\right),\quad x\in \K.
\]
Then $w_1Sw_1^{-1}=N_{\alpha_2}$ and hence $S\subset R_z$.
\end{itemize}

In conclusion, we observe that every  $R$-orbit in $Z$ contains an element $z$ with the following property: the stabilizer $R_z$ contains a subgroup $S_z$ such that the trivial representation of $S_z$ does not occur in \eqref{symp} as an irreducible subquotient, for all $k=0,1,2, \cdots$, and all but countably many $s\in \C$. Thus the trivial representation of $R_z$ does not occur  in \eqref{symp}  for all $k=0,1,2, \cdots$, and all but countably many $s\in \C$. 
This finishes the proof of the claim. 

By Theorem 1.15 and Example 1.16 of \cite{CS21}, for all but countably many $s\in \BC$ we have topological linear isomorphisms 
\[
\begin{aligned}
\RH^\CS_i(R; I(\sigma)\otimes\psi_s^{-1})&= \RH^\CS_i (R; \Gamma^\zeta(M, {\sf E})\otimes \psi_s^{-1})\\
&\cong  \RH^\CS_i (R; \Gamma^\zeta(U, \left.{\sf E}\right|_U)\otimes \psi_s^{-1})\\
& \cong \RH^\CS_i(R; \CS(R)\otimes\psi_s^{-1})\\
&\cong 
\left\{\begin{array}{ll} 
\BC, &  {\rm if \ } i=0;\\
\{0\}, &  {\rm if \ }i\neq 0.
\end{array}\right.
\end{aligned}
\]
\end{proof}

\section{Proof of Theorem \ref{thm2}}

This section is devoted to the proof of Theorem \ref{thm2}. We adopt the notations in previous sections. 

Let $f\in I(\sigma)$. We first assume that  $f|_{w_1R}\in \CS(w_1 R)$.  
Then $\La_s(f)$ converges absolutely and uniformly on compact subsets of $\BC$, and hence it is entire. It is clear that  for every $s\in \BC$ there exists such a function $f$ with $\La_s(f)\neq 0$.

Proposition \ref{prop3.1} implies that $f|_N\in \CS(N)$. We unfold $\oZ_s(f)$ as
\begin{eqnarray*}
\oZ_s(f)&=&\int_{\K^\times}\int_Nf(u[a])\psi_N^{-1}(u)|a|_\K^{s-\frac{n-1}{2}}\dd u\dd^\times a\\ 
&=&\int_{\K^\times}\int_N f([a]\cdot[a]^{-1}u[a])\psi_N^{-1}(u) |a|_\K^{s-\frac{n-1}{2}}\dd u\dd^\times a\\
&=&\int_{\K^\times}\int_N\si_1(a) f([a]^{-1}u[a])\psi_N^{-1}(u)|a|_\K^{s-(n-1)}\dd u\dd^\times a\\
&=&\int_{\K^\times}\int_{N'}\int_{\K}\si_1(a)f([a]^{-1}u'\cdot \mathrm u(x)[a])\psi_N^{-1}(u'\cdot \mathrm u(x)) |a|_\K^{s-(n-1)}\dd x \dd u'\dd^\times a.
\end{eqnarray*}
Note that $\psi_N$ restricted to $N'$ is invariant under the change of variable $u'\mapsto [a]u'[a]^{-1}$, and we have
\[
\big| \det({\rm Ad}([a])|_{{\rm Lie} \ N}) \big|_\K=|a|_\K^{n-1}.
\]
Hence
\[
\oZ_s(f)=\int_{\K^\times}\int_{N'}\int_{\K}\si_1(a)f(u'\cdot \mathrm u(x))\psi_N^{-1}(u')\psi^{-1}(ax) |a|_\K^{s}\dd x \dd u'\dd^\times a.
\]
For $u'\in N'$, define
\[
f_{u'}(x):=f(u'\cdot \mathrm u(x)),\quad x\in \K.
\]
Then $f_{u'} \in \CS(\K) $ and we have that
\begin{equation}\label{z}
\oZ_s(f)=\int_{\K^\times}\int_{N'} \CF_{\psi^{-1}}(f_{u'})(a) \sigma_1(a)|a|_\K^s \psi^{-1}_{N}(u')\dd u' \dd^\times a.
\end{equation}

Note that the function 
\[
 \left(N'\times \K\rightarrow \C,\quad (u',a)\mapsto  \CF_{\psi^{-1}}(f_{u'})(a)\right)\in \CS(N'\times \K).
\]
Thus by Tate's thesis, the integral \eqref{z} is absolutely convergent when $\mathrm{Re}(s)>-\nu_1$. 
Assume that this is the case. Then we can exchange the order of integration to obtain that
\begin{eqnarray*}
\oZ_s(f) & = & \int_{N'}\left( \int_{\K^\times} \CF_{\psi^{-1}}(f_{u'})(a) \sigma_1(a)|a|_\K^s \dd^\times a \right) \psi_N^{-1}(u')\dd u' \\
& = & \int_{N'} \oZ(s, \sigma_1,  \CF_{\psi^{-1}}(f_{u'}))  \psi_N^{-1}(u')\dd u'.
\end{eqnarray*}
By  \eqref{fe},  \eqref{eq2} and the Fourier inversion,
\begin{eqnarray*}
\oZ_s(f) & = & \gamma(s, \sigma_1, \psi)^{-1}   \int_{N'}  \oZ(1-s, \sigma_1^{-1}, f_{u'}) \psi_N^{-1}(u') \dd u' \\
& = & \gamma(s, \sigma_1, \psi)^{-1}  \int_{N'} \int_{\K^\times} f(u'\cdot \mathrm u(a))\sigma_1^{-1}(a)|a|_\K^{1-s} \psi_N^{-1}(u') \dd^\times a \dd u' \\
& = & \gamma(s, \sigma_1, \psi)^{-1} \int_{N'} \int_{\K^\times} f(\mathrm u(a)\cdot u')\sigma_1^{-1}(a)|a|_\K^{1-s} \psi_N^{-1}(u') \dd^\times a \dd u'  \\
& = &  \gamma(s, \sigma_1, \psi)^{-1} \cdot \La_s(f),
\end{eqnarray*} 
where in the second last equality we have made the change of variable $u'\mapsto {\rm Ad}(\mathrm u(a))(u')$ and used invariance property of $\psi_N$. 

Recall from Proposition \ref{prop1.1} that $\oZ_s(f)/ \oL(s, \sigma)$ is an entire function on $s\in \C$. Note that
\[
\frac{1}{\gamma(s, \sigma_1, \psi)\cdot \oL(s,\sigma)}=\frac{1}{\vep(s,\sigma_1,\psi)\cdot \oL(1-s, \sigma_1^{-1})\cdot \prod^n_{i=2}\oL(s, \sigma_i)},
\]
which is also an entire function on $s\in \C$.  
It follows from the uniqueness of  holomorphic continuation that
\begin{equation} \label{eq3}
\frac{\oZ_s(f)}{\oL(s, \sigma)}= \frac{\La_s(f)}{\gamma(s, \sigma_1, \psi)\cdot \oL(s,\sigma)}
\end{equation}
for any $s\in \BC$.

We now consider an arbitrary $f\in I(\sigma)$. By Theorem \ref{unique},    we have that  \eqref{eq3} holds for all but countably many $s\in \Omega_\sigma$, hence it holds for all $s\in \Omega_\sigma$ by continuity. This proves that
 \[
 \La_s(f)=\gamma(s, \sigma_1, \psi)\cdot\oZ_s(f)
 \]
 for  all $s\in \Omega_\sigma$ and $f\in I(\sigma)$.


\section*{Acknowledgements}

B. Sun was supported in part by National Key $\textrm{R}\,\&\,\textrm{D}$ Program of China (No. 2020YFA0712600) and  the National Natural Science Foundation of China (No. 11688101).

\end{document}